\documentclass[12pt]{amsart}
\usepackage[sc]{mathpazo}
\usepackage[centering,scale=0.8]{geometry}

\usepackage{amssymb,mathrsfs,graphicx,xcolor,microtype}
\usepackage[colorlinks,plainpages=false,
 linkcolor=blue!80!green,citecolor=green!60!blue,
 urlcolor=red!50!black]{hyperref}
\counterwithin{equation}{section}
\newtheorem{Thm}{Theorem}

\newtheorem{theorem}{Theorem}[section]
\newtheorem{lemma}[theorem]{Lemma}
\newtheorem{coro}[theorem]{Corollary}
\newtheorem{prop}[theorem]{Proposition}
\theoremstyle{definition}
\newtheorem{definition}[theorem]{Definition}
\theoremstyle{remark}
\newtheorem{example}[theorem]{Example}
\newtheorem{remark}[theorem]{Remark}
\newcommand{\Gr}{\operatorname{Gr}}
\newcommand{\Hom}{\operatorname{Hom}}

\newcommand{\CH}{\operatorname{CH}}
\def\aGr{\operatorname{\mathcal{G}\it r}}
\newcommand{\stacks}[1]{\href{https://stacks.math.columbia.edu/tag/#1}{Tag #1}}
\allowdisplaybreaks

\begin{document}
\title{A proof of the Grothendieck join--meet positivity}

\author{Peter L. Guo}
\address[Peter L. Guo]{Center for Combinatorics, LPMC, 
Nankai University, Tianjin 300071, P.R. China}
\email{lguo@nankai.edu.cn}

\author{Mingyang Kang}
\address[Mingyang Kang]{Center for Combinatorics, LPMC, 
Nankai University, Tianjin 300071, P.R. China}
\email{2120240005@mail.nankai.edu.cn}

\author{Jiaji Liu}
\address[Jiaji Liu]{Center for Combinatorics, LPMC, 
Nankai University, Tianjin 300071, P.R. China}
\email{1120240006@mail.nankai.edu.cn}

\author{Rui Xiong}
\address[Rui Xiong]{Department of Mathematics, University of Michigan, Ann Arbor, MI 48109}
\email{xiongrui@umich.edu}

\date{}
\maketitle

\begin{abstract}
We prove the (signed) Grothendieck  join--meet positivity 
$$
G_{\lambda\vee\mu}G_{\lambda\wedge\mu}
 -G_\lambda G_\mu
\in 
\sum_{\nu} \mathbb{N}\cdot (-1)^{|\nu|-|\lambda|-|\mu|}G_\nu
$$
for stable Grothendieck polynomials using geometric Satake equivalence, thereby  resolving a conjecture posed by Thomas and Yong. 
\end{abstract}

\setcounter{tocdepth}{1}
\tableofcontents

\section{Introduction}

For a partition $\lambda$, let $G_\lambda$ be the associated  stable Grothendieck polynomial which can be combinatorially defined via  set-valued tableaux of shape $\lambda$; see Buch  \cite[\S 3]{Buch}.
For two partitions $\lambda$ and $\mu$, their join $\lambda\vee \mu$ and meet $\lambda\wedge \mu$ are  defined respectively by  
$$
 (\lambda\vee\mu)_i=\max(\lambda_i,\mu_i),\qquad
 (\lambda\wedge\mu)_i=\min(\lambda_i,\mu_i).
$$
For example, for $\lambda=(3,1,1)$ and $\mu=(2,2,0)$, we have $\lambda\vee\mu=(3,2,1)$ and $\lambda\wedge\mu=(2,1,0)$. 

The main purpose of this paper is to prove the following (signed) Grothendieck positivity, which  was conjectured by 
Thomas and Yong  \cite[Conjecture 9.2]{TY}; 
see also Mihalcea \cite{Mil}.  

\begin{Thm}[\text{\cite[Conjecture 9.2]{TY}}]\label{thm:main}
For  partitions $\lambda$ and $\mu$,  
\begin{equation}\label{eq:main}
G_{\lambda\vee\mu}G_{\lambda\wedge\mu}
 -G_\lambda G_\mu
\in 
\sum_{\nu} \mathbb{N}\cdot (-1)^{|\nu|-|\lambda|-|\mu|}G_\nu.
\end{equation}
\end{Thm}

Here are some remarks.
\begin{itemize}
    \item The lowest degree component in  $G_\lambda$ is equal to the Schur function $s_\lambda$. So,  taking the lowest degree components on both sides of  \eqref{eq:main},  Theorem \ref{thm:main} reduces to  the celebrated     Lam--Postnikov--Pylyavskyy (LPP)   Schur positivity \cite[Theorem 5]{LPP}:
    \begin{equation}\label{schur-v}
s_{\lambda\vee\mu}s_{\lambda\wedge\mu}
 -s_\lambda s_\mu
\in 
\sum_{\nu} \mathbb{N}\cdot s_\nu.
\end{equation}

    \item Recently, significant progress toward Theorem \ref{thm:main} has   been made by  Chan,  Chen,   Pak  and   Soskin \cite[Theorem 1.8]{Chan} by developing the Ahlswede--Daykin--Schur inequality,  where a weaker version was proved, that is,  $G_{\lambda\vee\mu}G_{\lambda\wedge\mu}
 -G_\lambda G_\mu$ is (signed) Schur positive; see also \cite[Theorem 1.5]{GKL} for a refined version.
\end{itemize}

The proof of the LPP inequality   \eqref{schur-v} in  \cite{LPP} mainly relies on  Haiman's Schur positivity for  the  Kazhdan--Lusztig immanant \cite{Haiman} and   Rhoades--Skandera's result on Temperley--Lieb immanants \cite{RS}. 
In view of the arguments  in \cite{LPP},  the Schur positivity in \eqref{schur-v} can also  be deduced by  using \cite[Theorem 3]{Lam}, an identity of product of two Pl\"ucker coordinates and Temperley--Lieb immanants, together with Haiman's Schur positivity.  

The machinery behind the proof of the LPP Schur positivity  does not seem to extend to   Theorem \ref{thm:main}.
The key step in our proof is to give the identity a geometric representation-theoretic interpretation via the geometric Satake equivalence.

Geometrically, $G_\lambda$ represents the structure sheaf of the opposite Schubert variety $Y(\lambda)$ in the Grothendieck group of the Grassmannian $\Gr(k,n)$, provided that $\lambda$ lies  inside the rectangle $(n-k)^k$. 
The product $G_\lambda G_\mu$ will represent the structure sheaf of their transversal intersection. 

A direct attempt of proving Theorem \ref{thm:main} is by carefully degenerating  the intersection such that one is ``smaller'' than the other. 
But the intersection is usually hard to control. 
Instead, we will degenerate the product of two opposite Schubert varieties $Y(\lambda)\times Y(\mu)$ such that 
the degeneration of it will 
  be ``smaller'' than the degeneration of $Y(\lambda\vee\mu)\times Y(\lambda\wedge \mu)$.

Here is the degeneration principle used in the proof. For a
rank-$r$ vector bundle $\mathcal V$ on a projective complex scheme
$\mathfrak X$, define the \emph{Hirzebruch lambda class} 
$$
\mathsf\Lambda(y;\mathcal V)
:=\sum_{j=0}^r(-1)^j(1-y)^j[\wedge^j\mathcal V^\vee]\in K(\mathfrak{X})[y].$$
Let $x=(x_1,x_2,\ldots,x_k)$ be a family of variables. Write 
$$ \mathsf\Lambda(x;\mathcal V)
:=\prod_{i=1}^k\mathsf\Lambda(x_i;\mathcal V),\qquad 
 \chi(x;\mathcal V)
   :=\chi\bigl(\mathfrak X,\mathsf\Lambda(x;\mathcal V)\bigr)
       \in\mathbb Z[x_1,\ldots,x_k]^{S_k},
$$
where the Euler characteristic is taken coefficientwise. 

\begin{Thm}\label{thm:criterion}
Let $\mathcal Q$ be a globally generated rank-$r$ bundle on a
projective complex scheme $X$. Let
$\mathcal Y_\pm\subset U\times X$, where $U\subset\mathbb A^1$
is an open neighbourhood of $0$. Assume that 
\begin{enumerate}
 \item $\mathcal Y_\pm$ are integral, closed and flat over $U$,
 of relative dimension $d\le kr$;
 \item $\mathcal Y_\pm$ have rational singularities and integral
 geometric generic fibres;
 \item the degree-$(kr-d)$ component of
 $\chi(x;\mathcal Q|_{\mathcal Y_{\pm,0}})$ is nonzero; 
 \item $\mathcal Y_{-,0}\subseteq\mathcal Y_{+,0}$ as closed schemes.
\end{enumerate}
Then
\begin{equation}
\label{eq:criterion}
\chi(x;\mathcal Q|_{\mathcal Y_{+,0}})
-\chi(x;\mathcal Q|_{\mathcal Y_{-,0}})
\in\sum_{\nu\subseteq r^k}\mathbb{N}\cdot
(-1)^{|\nu|-(kr-d)}
G_\nu(x_1,\ldots,x_k).
\end{equation}
\end{Thm}

Now let us explain how to apply Theorem \ref{thm:criterion} to a proof of Theorem \ref{thm:main}. 
The degeneration family we used is from geometric Satake equivalence \cite{MV,BKK}. 

Write  $n=k+m$. We will work with the affine Grassmannian $\aGr$ of $GL_n$. 
The classical Grassmannian can be embedded into the affine Grassmannian, as an affine Schubert variety, 
\begin{equation}\label{eq:embGrtoaGr}
\Gr(k,n)\cong 
\aGr^{\varpi_m}
=\overline{\aGr^{\varpi_m}}
\stackrel{\subset}\longrightarrow \aGr.
\end{equation}
To define a family of schemes over $\mathbb{A}^1$, we still need the Beilinson--Drinfeld Grassmannian $\aGr_{\mathbb{A}}\to \mathbb{A}^1$ of $G$. 
Let $\mathbb{U}=\mathbb{A}^1\setminus \{0\}$. 
We have 
\begin{equation}\label{eq:BDGrtri}
\aGr_{\mathbb{U}}:= \mathbb{U}\times_{\mathbb{A}^1} \aGr_{\mathbb{A}}
\cong  \mathbb{U}\times \aGr\times \aGr. 
\end{equation}
Let $Y(\lambda)$ be the opposite Schubert variety of $\Gr(k,n)$. 
We define 
$$\mathcal{Y}_{\lambda,\mu} = \overline{\mathbb{U}\times Y(\lambda)\times Y(\mu)}
\subseteq \aGr_{\mathbb{A}},$$
where we used the identifications \eqref{eq:embGrtoaGr} and \eqref{eq:BDGrtri}. 
For readers familiar with geometric Satake equivalence, one might realize that $\mathcal{Y}_{\lambda,\mu}$ is a global version of Mirkovi\'c--Vilonen (MV) cycle used in \cite[Section 7.5]{BKK}. 

We will apply Theorem \ref{thm:criterion} to 
$$
\mathcal{Y}_+=\mathcal Y_{\lambda\vee\mu,\lambda\wedge\mu},
 \qquad \mathcal{Y}_-=\mathcal Y_{\lambda,\mu}.
$$
Note that the first three properties (1)--(3) in Theorem \ref{thm:criterion} are technical, and the most essential property is (4). 
Let us explain a weaker version of (4), namely, the cycle
$[\mathcal{Y}_{+,0}]-[\mathcal{Y}_{-,0}]$
is effective. 
By \cite[Section 7.5]{BKK}, the family 
$\mathcal{Y}_{\lambda,\mu}$ can be used to construct the natural map $L(\varpi_m)\otimes L(\varpi_m)\to L(2\varpi_m)$ geometrically.  
The class of the zero fiber is expanded positively in terms of MV cycles. 

It suffices to show 
\begin{equation}\label{eq:repofjointmeet}
\Delta_{C\vee D}\Delta_{C\wedge D}-\Delta_C\Delta_D \in \sum \mathbb{N}\{
\text{MV cycles of $L(2\varpi_m)$}\},
\end{equation}
where $\Delta_C$ is the MV cycle of $L(\varpi_m)$ indexed by an $m$-subset $C$ of $[n]$. It is known that \eqref{eq:repofjointmeet} is positive with respect to a combinatorial basis, called the \emph{Temperley--Lieb (TL) basis}; see \cite{Lam}. 
So the effectivity follows if we can prove the TL basis coincides with the MV cycles.

The paper is organized as follows. In Section \ref{sec:criterion-I}, we prove
Theorem \ref{thm:criterion}, a degeneration criterion for signed Grothendieck
positivity. In Section \ref{sec:main-I}, we construct the degeneration families
arising from geometric Satake and compare the Temperley--Lieb
and Mirkovi\'c--Vilonen bases to prove Theorem \ref{thm:main}.
Appendix \ref{sec:AppA} establishes the rational singularities of the families
used in the proof. Appendix \ref{sec:perfect} proves the uniqueness statement
for perfect bases needed for the comparison.

\subsection*{Acknowledgments}
The first author  was supported by the National Natural Science Foundation of China (No. 12371329) and the Fundamental Research Funds for the Central Universities (No. 63263094).

\subsection*{AI disclosure}
This project began in early June 2026. Over the following months, the authors engaged in extensive discussions with ChatGPT, which, in late September 2026, provided the main strategy for the proof of Theorem A. The authors have since substantially revised, expanded, and simplified the resulting  material and take full responsibility for the final content of this manuscript.

\section{Proof of Theorem B}\label{sec:criterion-I}

\subsection{A geometric meaning of $\chi(x;\mathcal V)$}
Let $\mathcal V$ be a rank-$r$ bundle on a quasi-projective
scheme $\mathfrak X$, generated by a finite-dimensional space $H$, i.e. 
$\mathcal O_{\mathfrak X}\otimes H\twoheadrightarrow\mathcal V$.
Fix a generic linear map
$$
 \phi:\mathbb C^{k+r}\longrightarrow H.
$$
On $\mathfrak X\times\Gr(k,k+r)$, consider
$$
 s:\mathcal S\hookrightarrow\mathcal O^{k+r}
    \xrightarrow{\phi}\mathcal O\otimes H
    \twoheadrightarrow\mathcal V,
$$
where $\mathcal S$ is the tautological subbundle. Define
\begin{equation}
\label{eq:defofincidence}
\begin{aligned}
 \widetilde{\mathfrak X}
  &:=Z(s)
   =\{(z,V):\phi(V)\text{ vanishes at }z\},\\
 \pi&:\widetilde{\mathfrak X}\longrightarrow\Gr(k,k+r)
       \quad\text{(second projection).}
\end{aligned}
\end{equation}
Bundles on the two factors are pulled back without changing notation.

\begin{example}\label{eg:tildeX=X}
Take $H=H^0(\mathfrak{X},\mathcal{V})$.  When $k+r\geq \dim H$, the generic choice $\phi$ is surjective. 
In this case, $\tilde{\mathfrak{X}}\cong \mathfrak{X}$ via the first projection, and the second projection 
$\mathfrak{X}\cong \tilde{\mathfrak{X}}\to \Gr(k,k+r)$ coincides with the morphism defined by the universal property of $\Gr(k,k+r)$. Here, we can view $\mathcal{V}$ as a quotient bundle of the trivial bundle $\mathcal{O}^{k+r}$. 
\end{example}

\begin{example}
Under the hypotheses of Theorem \ref{thm:eq:geoofchi} below,
$\widetilde{\mathfrak X}$ is the closure of the graph of the rational map
$$
 \mathfrak X\dashrightarrow\Gr(k,k+r),\qquad
 z\longmapsto\ker\bigl[
 \mathbb C^{k+r}\xrightarrow{\phi}H
       \xrightarrow{\operatorname{ev}_z}\mathcal V_z\bigr].
$$
This map is defined on the open set where
$\operatorname{ev}_z\circ\phi$ is surjective.
\end{example}

\begin{theorem}\label{thm:eq:geoofchi}
Assume that $\mathfrak X$ is a variety with rational singularities. For generic
$\phi$, so is $\widetilde{\mathfrak X}$, and the first
projection is birational. If $\mathfrak X$ is projective, then
\begin{equation}\label{eq:geoofchi}
 \pi_*[\mathcal O_{\widetilde{\mathfrak X}}]
 =\chi(1-[\mathcal S_1],\ldots,1-[\mathcal S_k];\mathcal V)
\end{equation}
in $K(\Gr(k,k+r))$, where $\mathcal S_i$ are the formal
$K$-theoretic Chern roots of $\mathcal S$.
\end{theorem}
\begin{proof}
Allow $\phi$ to vary. The universal incidence $\mathcal X$ is
the total space of the vector bundle
$$
 \ker\left[
 \Hom(\mathbb C^{k+r},H)\otimes\mathcal O
 \twoheadrightarrow\mathcal S^\vee\otimes\mathcal V
 \right]
$$
on $\mathfrak X\times\Gr(k,k+r)$. Thus $\mathcal X$ has
rational singularities. With $d=\dim\mathfrak X$,
$$
 \dim\mathcal X=d+\dim\Hom(\mathbb C^{k+r},H).
$$
Its projection to $\Hom(\mathbb C^{k+r},H)$ is surjective:
every map $\mathbb C^{k+r}\to\mathcal V_z$ has a kernel
containing a $k$-plane. Hence a general fibre
$\widetilde{\mathfrak X}$ has pure dimension $d$ and rational
singularities \cite[Lemma 3]{BrionPos}.

On the universal incidence scheme, the failure of surjectivity is 
$$
\det\bigl[\mathbb C^{k+r}/\mathcal S
\longrightarrow\mathcal V\bigr]=0.
$$
Its fibre over a general $\phi$ has dimension at most $d-1$. 
Outside this locus, the first projection identifies
$\widetilde{\mathfrak X}$ with the graph of
$$
z\longmapsto
\ker\bigl[\mathbb C^{k+r}\longrightarrow\mathcal V_z\bigr]
$$
over the open subset of $\mathfrak X$ where the map is surjective.
Since $\widetilde{\mathfrak X}$ is pure of dimension $d$,
this graph is nonempty and dense. The graph is integral,
and $\widetilde{\mathfrak X}$ is reduced because it has
rational singularities. Consequently,
$\widetilde{\mathfrak X}$ is integral, and its first projection
to $\mathfrak X$ is birational.

Finally, rational singularities imply that $\mathfrak X$ is
Cohen--Macaulay \cite[notation and conventions, p. 3]{BrionPos}. Since
$$
 \operatorname{codim}\widetilde{\mathfrak X}=kr
       =\operatorname{rank}(\mathcal S^\vee\otimes\mathcal V),
$$
$s$ is regular \cite[\stacks{02JN}]{Stacks}. Its Koszul
resolution gives, in $K(\mathfrak X\times\Gr(k,k+r))$,
\begin{align*}
 [\mathcal O_{\widetilde{\mathfrak X}}]
  &=\sum_{j\geq 0}(-1)^j[\wedge^j(\mathcal S\otimes\mathcal V^\vee)]
  =\mathsf\Lambda(1-[\mathcal S_1],\ldots,
                    1-[\mathcal S_k];\mathcal V).
\end{align*}
Applying the second projection's pushforward proves
\eqref{eq:geoofchi}.
\end{proof}

We use the finite-rectangle identification
\begin{align}\label{eq:rectangle}
 &\{P\in\mathbb Z[x_1,\ldots,x_k]^{S_k}:\deg_{x_i}P\le r\}
       \xrightarrow{\ \sim\ }K(\Gr(k,k+r)),\notag\\
 &\hspace{22mm}P\longmapsto P(1-[\mathcal S_1],\ldots,1-[\mathcal S_k]).
\end{align}
This is an isomorphism of abelian groups. For $\lambda\subseteq r^k$,
$$
 \deg_{x_i}G_\lambda\le\lambda_1\le r,\qquad
 G_\lambda=s_\lambda+\text{higher-degree terms}.
$$
Thus these polynomials form a basis and map to the Schubert basis
$[\mathcal O_{Y(\lambda)}]$; see
\cite[Theorems 3.1, 8.1 and (8.1)]{Buch}.

\begin{coro}\label{cor:schubert-euler}
For $\lambda\subseteq(n-k)^k$, let $\mathcal Q$ be the universal
quotient bundle on $\Gr(k,n)$. Then
$$
 \chi(x;\mathcal Q|_{Y(\lambda)})=G_\lambda(x_1,\ldots,x_k).
$$
\end{coro}
\begin{proof}
Take $H=\mathbb C^n$ and $\phi=\mathrm{id}$. The incidence is
the graph of $Y(\lambda)\hookrightarrow\Gr(k,n)$.
The same Koszul calculation gives
$$
 \chi(1-[\mathcal S_1],\ldots,1-[\mathcal S_k];
              \mathcal Q|_{Y(\lambda)})
       =[\mathcal O_{Y(\lambda)}].
$$
Both polynomials have degree at most $n-k$ in each variable.
Apply \eqref{eq:rectangle} with $r=n-k$.
\end{proof}

\begin{lemma}\label{lem:leading-degree}
Let $\mathfrak X$ be a projective variety with rational singularities and
$d=\dim\mathfrak X\le kr$. If the degree-$(kr-d)$ component of
$\chi(x;\mathcal V)$ is nonzero, then $\pi$ is generically finite
onto its image.
\end{lemma}
\begin{proof}
Pick a resolution of $\widetilde{\mathfrak X}$. 
Rational singularities and
Grothendieck--Riemann--Roch \cite[Theorem 15.2]{Fulton} give that 
$\pi_*[\widetilde{\mathfrak{X}}]\neq 0$ in the Chow group
$\CH^{kr-d}(\Gr(k,k+r))$.
Thus $\pi:\widetilde{\mathfrak{X}}\to\pi(\widetilde{\mathfrak{X}})$ is generically finite.
\end{proof}

\subsection{Interpretation of the coefficients}
Let $X(\nu)\subset\Gr(k,k+r)$ be the Schubert variety dual to
$Y(\nu)$, and let $\mathcal D=\det\mathcal S^\vee$ be the determinant bundle. Our indexing is
$$
 \operatorname{codim}Y(\nu)=|\nu|,\qquad \dim X(\nu)=|\nu|.
$$
The Euler-dual basis formula is
\begin{equation}\label{eq:duality}
 \chi\bigl(\Gr(k,k+r),[\mathcal O_{Y(\lambda)}]
       [\mathcal O_{X(\nu)}][\mathcal D^{-1}]\bigr)
       =\delta_{\lambda,\nu};
\end{equation}
see \cite[\S8, paragraph preceding (8.2)]{Buch}.

\begin{lemma}[Kleiman--Bertini]\label{lem:transversality}
Let $f:Y\to\Gr(k,k+r)$ be a morphism from a variety with rational
singularities. For generic $g\in GL_{k+r}$,
$$
 W=f^{-1}(gX(\nu))
$$
is empty or has rational singularities, of pure dimension
$\dim Y+|\nu|-kr$. Moreover,
$$
 \operatorname{Tor}^{\Gr(k,k+r)}_i
       (\mathcal O_Y,\mathcal O_{gX(\nu)})=0\qquad(i>0).
$$
\end{lemma}
\begin{proof}
The action of $GL_{k+r}$ is transitive, and its orbit map has
Zariski local sections; see \cite[Definition 2.1, Proposition 2.2]{BCMP}.
Since $X(\nu)$ has rational singularities, apply
\cite[Theorem 2.5]{BCMP} for rational singularities,
\cite[Theorem 2(i)]{Kleiman} for dimension, and
\cite[Theorem 2.1]{Sierra} for Tor vanishing.
The dimension and Tor assertions also apply to pure-dimensional
schemes without rational singularities.
\end{proof}

Assume $\mathcal Y\subset U\times X$ satisfies (1)--(3) of
Theorem \ref{thm:criterion}. Using a fixed generating space
$H=H^0(X,\mathcal Q)$, form
$$
 \pi:\widetilde{\mathcal Y}\longrightarrow\Gr(k,k+r)
$$
using \eqref{eq:defofincidence}. 
The total incidence is integral with rational singularities.
Since it dominates the smooth curve $U$, it is flat. Thus
$$
 \widetilde{\mathcal Y_0}=(\widetilde{\mathcal Y})_0,
 \qquad\dim\widetilde{\mathcal Y}_0=d.
$$
The Koszul resolution on the total space remains exact on the
fibre: its terms and its resolved sheaf are flat over $U$.
Consequently \eqref{eq:geoofchi} holds on the special fibre,
without assuming that this fibre has rational singularities.
For generic $g\in GL_{k+r}$, write
$$
 \chi(x;\mathcal Q|_{\mathcal Y_0})=\sum_\nu b_\nu G_\nu(x),
 \qquad
 \mathcal W=\pi^{-1}(gX(\nu)),\quad W_0=\mathcal W_0.
$$

\begin{theorem}\label{thm:coefficients}
We have 
$b_\nu=\chi(W_0,\pi^*\mathcal D^{-1})$.
\end{theorem}
\begin{proof}
Apply \eqref{eq:duality}, the projection formula and Tor vanishing:
\begin{align*}
 b_\nu
 &=\chi\bigl(\Gr(k,k+r),\pi_{0*}[\mathcal O_{\widetilde{\mathcal Y}_0}]
             [\mathcal O_{gX(\nu)}][\mathcal D^{-1}]\bigr)\\
 &=\chi\bigl(\widetilde{\mathcal Y}_0,
              [\mathcal O_{W_0}]\,\pi^*[\mathcal D^{-1}]\bigr)
  =\chi(W_0,\pi^*\mathcal D^{-1}).\qedhere
\end{align*}
\end{proof}

\subsection{A vanishing Theorem}
For a line bundle $\mathcal{L}$ on a projective variety $W$, let us recall the following concepts. 
The line bundle $\mathcal{L}$ is \emph{nef} if
 $\deg(\mathcal{L}|_C)\ge0$
 for every integral curve $C\subset W$. 
More generally, for a morphism $f:W\to U$, the line bundle $\mathcal{L}$ is \emph{$f$-nef} if it has nonnegative degree on curves in fibres. 
The line bundle $\mathcal{L}$ is \emph{big} if $$\displaystyle\limsup_{m\to\infty}\frac{\dim 
H^0(W,\mathcal{L}^m)}{m^{\dim W}}>0.$$
The implications used below are
 \begin{align*}
  \mathcal{L}\text{ globally generated}&\ \Longrightarrow\ \mathcal{L}\text{ nef},\\
  \mathcal{L}\text{ ample}&\ \Longrightarrow\ L\text{ nef and big},\\
  f\text{ generically finite},\ A\text{ ample}
    &\ \Longrightarrow\ f^*A\text{ nef and big}.
\end{align*}

\begin{lemma}\label{lem:relative-vanishing}
Let $f:W\to U$ be a projective flat family of pure relative
dimension $e$ over a smooth complex curve. Assume $W$ has rational
singularities. If $\mathcal{L}$ is $f$-nef and big on every component of the
geometric generic fibre, then
$$
 H^j(W_0,\mathcal{L}_0^{-1})=0\qquad(0\le j<e),
$$
where $W_0=f^{-1}(0)$ and $\mathcal{L}_0=\mathcal{L}|_{W_0}$.
\end{lemma}
\begin{proof}
Take a projective resolution $\rho:\widetilde W\to W$, working
componentwise. Rational singularities, Grauert--Riemenschneider
vanishing and relative Kawamata--Viehweg vanishing give
\begin{align*}
 R\rho_*\omega_{\widetilde W}&=\omega_W,\\
 R^i(f\rho)_*(\omega_{\widetilde W}\otimes\rho^*\mathcal{L})&=0\qquad(i>0).
\end{align*}
For the canonical-sheaf identification, see
\cite[Notation and conventions, p. 3]{BrionPos}; for vanishing,
see \cite[Theorem 1.2]{FujinoNote}. Here $\rho^*L$ is relatively
nef and big on the generic fibre. Projection formula gives
\begin{equation}\label{eq:omega-direct}
 R^if_*(\omega_W\otimes \mathcal{L})=0\qquad(i>0).
\end{equation}

Shrink $U$ around $0$ and choose a parameter $s$ cutting out $0$.
Since $W$ is Cohen--Macaulay and flat over $U$, Cartier adjunction
gives
$$
 0\longrightarrow\omega_W\otimes \mathcal{L}
 \stackrel{s}\longrightarrow 
 \omega_W\otimes \mathcal{L}
 \longrightarrow\omega_{W_0}\otimes \mathcal{L}_0\longrightarrow0,
$$
where the last sheaf is extended by zero to $W$;
see \cite[\stacks{0A7N}, \stacks{0AWS}]{Stacks}.
Applying $Rf_*$ and \eqref{eq:omega-direct}, we obtain
$$
 H^i(W_0,\omega_{W_0}\otimes L_0)=0\qquad(i>0).
$$
The Cartier fibre $W_0$ is Cohen--Macaulay of dimension $e$.
Serre duality \cite[\stacks{0FVZ}]{Stacks} now gives
$$
 H^j(W_0,\mathcal{L}_0^{-1})^*
 \simeq H^{e-j}(W_0,\omega_{W_0}\otimes \mathcal{L}_0)=0
 \qquad(j<e).\qedhere
$$
\end{proof}

\subsection{Proof of Theorem \ref{thm:criterion}}
Assume $\mathcal Y_+$ and $\mathcal{Y}_-$ satisfy (1)--(4) of
Theorem \ref{thm:criterion}, and write
$$
 \chi(x;\mathcal Q|_{\mathcal Y_{+,0}})
 -\chi(x;\mathcal Q|_{\mathcal Y_{-,0}})
       =\sum_\nu b_\nu G_\nu(x).
$$
We will prove $(-1)^{|\nu|-(kr-d)}b_\nu\ge0$.

\begin{proof}[Proof of Theorem \ref{thm:criterion}]
Take $H=H^0(X,\mathcal Q)$ and use the same generic
$\phi:\mathbb C^{k+r}\to H$ for both families. As in
Section \ref{sec:criterion-I}, the incidence families are flat
with rational singularities, and we denote the natural morphism by
$$
 \pi_\pm:\widetilde{\mathcal Y}_\pm\longrightarrow\Gr(k,k+r).$$
It follows from the construction of incidence family that 
$\widetilde{\mathcal Y}_{-,0}\subseteq\widetilde{\mathcal Y}_{+,0}$. 
Their geometric generic fibres have rational singularities by flat base change. 
They are
integral because they are birational to those of $\mathcal Y_\pm$.
By condition (3) and Lemma \ref{lem:leading-degree},
$$
 \pi_{\pm,\eta}\text{ is generically finite onto its image}.
$$

Fix a partition $\nu\subseteq r^k$ and choose the same generic $g$ for both signs. Put
$$
 e=|\nu|-(kr-d),\qquad
 \mathcal W_\pm=\pi_\pm^{-1}(gX(\nu)),\qquad
 W_{\pm,0}=(\mathcal W_\pm)_0.
$$
Generic transversality gives rational singularities of
$\mathcal W_\pm$ and, when nonempty,
$$
 \dim\mathcal W_\pm=e+1,\qquad\dim W_{\pm,0}=e.
$$
If $e<0$, both special tests are empty and $b_\nu=0$.
For $e\ge0$, no component of $\mathcal W_\pm$ lies over $0$.
Remove the other vertical components by shrinking $U$.
The remaining normal components dominate $U$, hence are flat
of pure relative dimension $e$.

Generic transversality proves 
$$\pi_{\pm,\eta}|_{W_{\pm,\eta}}
\text{ is generically finite onto its image}.$$
So $\pi_\pm^*\mathcal D$ is nef and big on every geometric generic-fibre component. 
It is also relatively nef, since $\mathcal D$ is globally generated.

Lemma \ref{lem:relative-vanishing} and dimension vanishing give
\begin{equation}\label{eq:concentration}
 H^j(W_{\pm,0},\pi_\pm^*\mathcal D^{-1})=0\qquad(j\ne e).
\end{equation}
Condition (4) gives $W_{-,0}\subseteq W_{+,0}$. Let $\mathcal I$
be its ideal sheaf. 
Consider the short exact sequence 
on $W_{+,0}$
$$
 0\longrightarrow\mathcal I\otimes\pi_+^*\mathcal D^{-1}
 \longrightarrow\pi_+^*\mathcal D^{-1}
 \longrightarrow\pi_-^*\mathcal D^{-1}\longrightarrow0
$$
where the last sheaf is extended by zero. Since
$$
 H^{e-1}(W_{-,0},\pi_-^*\mathcal D^{-1})=0,\qquad
 H^{e+1}(W_{+,0},\mathcal I\otimes\pi_+^*\mathcal D^{-1})=0,
$$
we obtain a short exact sequence
$$
0\longrightarrow
 H^e(W_{+,0},\mathcal I\otimes\pi_+^*\mathcal D^{-1})
 \longrightarrow H^e(W_{+,0},\pi_+^*\mathcal D^{-1})
 \longrightarrow H^e(W_{-,0},\pi_-^*\mathcal D^{-1})
 \longrightarrow0.$$
Finally, Theorem \ref{thm:coefficients}
and \eqref{eq:concentration} give
\begin{align*}
 b_\nu
  &=\chi(W_{+,0},\pi_+^*\mathcal D^{-1})
           -\chi(W_{-,0},\pi_-^*\mathcal D^{-1}),\\
 (-1)^e b_\nu
  &=\dim H^e(W_{+,0},\mathcal I\otimes\pi_+^*\mathcal D^{-1})
       \ge0.
\end{align*}
This proves Theorem \ref{thm:criterion}.
\end{proof}

\section{Proof of Theorem A}\label{sec:main-I}

In this section, we denote $n=k+m$ and $G=GL_n$. 
For a partition $\lambda\subseteq m^k$, set
\begin{equation}\label{eq:column-subset}
 C_\lambda=\{k+j-\lambda'_j:1\le j\le m\}\subseteq[n],
\end{equation}
where $\lambda'$ is the conjugate partition. 
The map $\lambda\mapsto C_\lambda$ defines a bijection between partitions inside the rectangle $m^k$ and $m$-subsets of $[n]$. 
Let us denote by $[\lambda]\in \mathbb{Z}^n$ the indicator vector of $C_\lambda$.  
Note that $Y(\lambda)$ is the $B^-$-orbit closure of $\operatorname{span}(e_i:i\notin C_\lambda)$.

\subsection{Identifying family}
\label{sec:satake}
Let us recall some background of Geometric Satake equivalence \cite{LusztigSatake,GinzburgSatake,MV}. 
The affine Grassmannian of $G=GL_n$ is
$$
 \aGr
 =\{\text{$\mathbb{C}[[t]]$-lattices }L\subset \mathbb{C}(\!(t)\!)^{\oplus n}\}. 
$$
Let $L_0=\mathbb{C}[[t]]^{\oplus n}\in \aGr$ be the standard lattice. 
For a dominant coweight
$\theta=(\theta_1\ge\cdots\ge\theta_n)\in\mathbb Z^n$, the spherical Schubert cell is 
$\aGr^\theta=G[[t]]t^\theta L_0$ where 
$t^\theta=\operatorname{diag}(t^{\theta_1},\ldots,t^{\theta_n})\in G(\!(t)\!):= G(\mathbb{C}(\!(t)\!))$.
Let $N^-\subset G$ be the lower unitriangular subgroup.
For $\gamma\in\mathbb Z^n$, denote the semi-infinite orbit by 
$S^-_\gamma=N^-(\!(t)\!)t^\gamma L_0\subset \aGr$. 

Geometric Satake equivalence \cite[Theorem 7.3]{MV} is a tensor equivalence 
\begin{equation}\label{eq:satake}
\operatorname{Perv}_{G[[t]]}(\aGr,\mathbb C)
\simeq\operatorname{Rep}_{\mathbb C}(G^\vee).
\end{equation}
The left-hand category consists of $G[[t]]$-equivariant perverse sheaves over $\aGr$, and the right-hand category is the category of finite-dimensional
representations of $G^\vee=GL_n=G$.

The intersection-cohomology
complex $\mathbf{IC}_\theta$ of $\overline{\aGr^\theta}$ corresponds
to the irreducible representation $L(\theta)$ of highest weight
$\theta$. 
Let $\gamma$ be a weight of $L(\theta)$, 
Mirkovi\'c and Vilonen \cite{MV} proved a geometric realization of the weight space 
$$
L(\theta)_\gamma\simeq 
H_{\mathsf{top}}(\aGr^\theta\cap S^-_\gamma,\mathbb C)\simeq 
H_{\mathsf{top}}(\overline{\aGr^\theta}\cap\overline{S^-_\gamma},\mathbb C).
$$
Here all homology groups are Borel--Moore homology, and 
$\aGr^\theta\cap S^-_\gamma$ is 
pure of dimension $\langle \rho,\theta-\gamma\rangle$ and ``${\sf top}$'' means the top degree, i.e. $2\langle \rho,\theta-\gamma\rangle$. 
In particular, the fundamental classes of irreducible components of 
$\overline{\aGr^\theta}\cap \overline{S^-_\gamma}$
form a basis of $L(\theta)_\gamma$. 
They are called the \emph{MV cycles}. 

\begin{example}
For $\varpi_m=(1^m,0^{k})$, we can identify 
$$
 \aGr^{\varpi_m}=\overline{\aGr^{\varpi_m}}
 =\{L\in \aGr: tL_0\subseteq L\subseteq L_0,\, \dim L_0/L=m\}. 
$$
Thus 
$$
 \aGr^{\varpi_m}=\overline{\aGr^{\varpi_m}}
 \stackrel{\sim}\longrightarrow\Gr(k,n),\qquad
 L\longmapsto L/tL_0
 \leq L_0/tL_0\simeq \mathbb{C}^n.
$$
For a partition $\lambda\subseteq m^k$, this map identifies
$$
\overline{\aGr^{\omega_m}}\cap \overline{S_-^{[\lambda]}}
 \simeq Y(\lambda)\subseteq \Gr(k,n). 
$$
\end{example}

\begin{example}
We have the following identification
$$
\begin{aligned}
{\aGr^{2\varpi_m}}
& =\{L\in \aGr: t^2L_0\subseteq L\subseteq L_0,\, L_0/L\simeq (\mathbb{C}[[t]]/(t^2))^{\oplus m}\}\\
\overline{\aGr^{2\varpi_m}}
&=\{L\in \aGr: t^2L_0\subseteq L\subseteq L_0,\dim L_0/L=2m\}.
\end{aligned}$$
\end{example}

For our purpose, we still need the \emph{Beilinson--Drinfeld Grassmannian}
$\aGr_{\mathbb{A}}\to\mathbb{A}$. 
Its fibre parametrizes rank-$n$ vector bundles
on $\mathbb{A}^1$ with a trivialization away from $\{0,s\}$
\cite[Definition 3.1.1]{Zhu}. 
Let $\mathbb{U}=\mathbb{A}^1\setminus \{0\}$. 
We have 
\begin{equation}\label{eq:BDGrtri2}
\aGr_{\mathbb{U}}:= \aGr_{\mathbb{A}}\times_{\mathbb{A}^1} \mathbb{U}
\cong  \mathbb{U}\times \aGr\times \aGr. 
\end{equation}
Denote by 
$\mathcal{X}=\overline{\aGr_{\mathbb{A}}^{\varpi_m,\varpi_m}}$ 
the closure of
$\mathbb{U}\times\aGr^{\varpi_m}\times\aGr^{\varpi_m}$
in $\aGr_{\mathbb{A}}$; see
\cite[Definition 3.9]{BDK}.
Its special fibre is isomorphic to $\overline{\aGr^{2\varpi_m}}$
as a scheme \cite[Proposition 3.1.14]{Zhu}.

Let us give an explicit description of $\mathcal{X}$. 
Let 
$$\mathcal{X}'=
\overline{\{V\oplus_s W: s\in \mathbb{U}, V\in \Gr(k,n),W\in \Gr(k,n)\}}\subseteq \mathbb{A}^1\times \Gr(2k,2n)$$
where 
$V\oplus_s W
 =\{(u,v)\in\mathbb{C}^{2n}:u\in V,\ u+sv\in W\}$.

\begin{theorem}\label{thm:bd-model}
We have an isomorphism over $\mathbb{A}^1$
\begin{equation}\label{eq:bd-identification}
\mathcal{X}:=\overline{\aGr_{\mathbb{A}}^{\varpi_m,\varpi_m}}
\stackrel{\sim}\longrightarrow  \mathcal{X}'. 
\end{equation}
\end{theorem}
\begin{proof}
By the proof of \cite[Theorem 3.1.3]{Zhu}, we have a closed embedding
$$\overline{\aGr_{\mathbb{A}}^{\varpi_m,\varpi_m}}
\subseteq \{(s,L'): t(t-s)L_0'\subseteq L'\subseteq L_0', \,\dim L'/t(t-s)L_0'=2k\},$$
where $L'$ is a $\mathbb{C}[t]$-submodule, and $L_0'=\mathbb{C}[t]^{\oplus n}$. 

Note that 
$\mathbb{C}^{n}\oplus \mathbb{C}^n\cong L_0'/t(t-s)L_0$ by $(x,y)\mapsto x+ty$. 
We thus get a closed embedding 
$$\overline{\aGr_{\mathbb{A}}^{\varpi_m,\varpi_m}}\longrightarrow \mathbb{A}^1\times \Gr(2k,2n),\qquad 
L'\longmapsto L'/t(t-s)L_0\leq L_0/t(t-s)L_0\cong \mathbb{C}^n\oplus \mathbb{C}^n. $$
For $s\neq 0$, we have a natural isomorphism of $\mathbb{C}[t]$-modules
$$L'_0/t(t-s)L_0' \simeq \mathbb{C}_{t=0}^{n}\oplus \mathbb{C}_{t=s}^{n}
\simeq 
\mathbb{C}[[t]]^{\oplus n}/t\mathbb{C}[[t]]^{\oplus n}\oplus 
\mathbb{C}[[t']]^{\oplus n}/t'\mathbb{C}[[t']]^{\oplus n}$$
where $t'=t-s$. This defines the embedding $\aGr^{\varpi_m}\times \aGr^{\varpi_m}\to \overline{\aGr}_{\mathbb{A}}^{\varpi_m,\varpi_m}$. 
We leave it to the reader to check this is an isomorphism over $\mathbb{U}$. 
By taking closure, we get the isomorphism. 
\end{proof}

\begin{prop}\label{prop:fusion-rational}
The variety $\mathcal Y_{\lambda,\mu}$ has rational singularities.
Every fibre is Cohen--Macaulay and pure of dimension
$2km-|\lambda|-|\mu|$.
\end{prop}
\begin{proof}
The proof that $\mathcal{Y}_{\lambda,\mu}$ has rational singularities is technical, we give the proof in the appendix, see Theorem \ref{prop:RSofYlambdaYmu}. 
Since $\mathcal{Y}_{\lambda,\mu}$ 
is integral and dominates $\mathbb A^1$,  it is flat over $\mathbb{A}^1$. 
Rational singularities imply
Cohen--Macaulayness, and $s-s_0$ is a non-zero-divisor for every
$s_0\in\mathbb C$. Hence each fibre is Cohen--Macaulay, of pure
dimension
$$
 \dim\mathcal Y_{\lambda,\mu}-1
 =\dim Y(\lambda)+\dim Y(\mu)
 =2km-|\lambda|-|\mu|.\qedhere
$$
\end{proof}

\subsection{Properties (1)--(3)}
Now let us set 
$$r=2m,\qquad U=\mathbb{A}^1,\qquad X=\Gr(2k,2n)$$
and $\mathcal{Q}$ the universal quotient bundle. 
By Theorem \ref{thm:bd-model}, we may realize
$$\mathcal{Y}_{\lambda,\mu}\subseteq \mathcal{X}\simeq \mathcal{X}'\subseteq U\times X.$$

Now we can prove properties (1)--(3). 

\begin{prop}\label{step:(1)--(3)}
We have 
$
 \chi(x;\mathcal Q|_{\mathcal{Y}_{\lambda,\mu,0}})
     =G_\lambda(x)G_\mu(x)$.
Moreover, the variety $\mathcal{Y}_{\lambda,\mu}$ satisfies the conditions (1)--(3) in Theorem \ref{thm:criterion}. 
\end{prop}
\begin{proof}
Since $\mathcal Y_{\lambda,\mu}|_{\mathbb{U}}
 \simeq\mathbb \mathbb{U}\times Y(\lambda)\times Y(\mu)$, property (1) is obvious.
Proposition \ref{prop:fusion-rational}
gives property (2). 

On the general fiber at $s$, the quotient
bundle is the direct sum of the two quotient bundles of the two factors. Therefore
$$ \chi(x;\mathcal Q|_{\mathcal Y_{\lambda,\mu,s}})
=\chi(x;\mathcal Q|_{Y(\lambda)})
           \chi(x;\mathcal Q|_{Y(\mu)})
    =G_\lambda(x)G_\mu(x).
$$
Here we used the general fact of Hirzebruch lambda class: 
$$\mathsf\Lambda(y;\mathcal V\oplus\mathcal U)
=\mathsf\Lambda(y;\mathcal V)\mathsf\Lambda(y;\mathcal U).$$
Flat invariance of Euler characteristic gives the same identity at $s=0$
\cite[III, Theorem 9.9]{Hartshorne}.

Finally, put $d=2km-|\lambda|-|\mu|$. By above, the $2km-d$ degree component of 
$\chi(x;\mathcal Q|_{\mathcal Y_{\lambda,\mu,\eta}})=G_\lambda(x)G_\mu(x)$ is $s_\lambda(x)s_\mu(x)\ne0$. 
The last inequality holds
since $\lambda,\mu$ have at most $k$ rows. This proves (3). 
\end{proof}

\subsection{The Temperley--Lieb and MV bases}
Let $L(\theta)$ be the finite-dimensional representation of $GL_n$ of highest weight $\theta\in \mathbb{Z}^n$. 
Recall that $L(\varpi_m)$ can be realized as the span of minors over $m\times n$ matrices; and 
$L(2\varpi_m)$ can be realized as the span of products of two minors (i.e. with Pl\"ucker relation). 

For $C$ an $m$-subset of $[n]$, let 
$\Delta_C$ be the corresponding minor. 
Let us call $\{\Delta_C\}$ the minor basis. 
The Temperley--Lieb basis $\{F_{\tau,D}\}_{(\tau,D)}$ is a basis of $L(2\varpi_m)$ indexed by pairs $(\tau,D)$, where $\tau$ is a perfect noncrossing matching over a subset $A\subset [n]$, and $D\subset [n]\setminus A$ is disjoint from $A$ such that $2|D|+|A|=2m$; see   Lam \cite[\S 3]{Lam}. 
Following  \cite[Theorem 3.1]{Lam}, for two $m$-subsets $C,C'\subset [n]$, 
\begin{equation}\label{eq:minortoTL}
\Delta_C\,\Delta_{C'}
 =\sum_{(\tau,D)\in\mathcal C(C,C')}F_{\tau,D},
\end{equation} 
where $\mathcal{C}(C,C')$ is the set of pairs $(\tau,D)$ such that $\tau$ matches $C\setminus C'$ with $C'\setminus C$, and $D=C\cap C'$.

\begin{lemma}\label{lem:tlmv}
The minor basis equals the MV cycles of $L(\omega_m)$. 
The Temperley--Lieb basis equals the MV cycles of $L(2\omega_m)$. 
\end{lemma}
\begin{proof}
The Temperley--Lieb basis is the upper global basis
\cite[Theorem 2.1(2) and Proposition 3.2]{Lam}, hence is perfect.
The MV basis is perfect by \cite[Theorem 5.2]{BKK}.
Since
$$
 \langle2\omega_m,\beta^\vee\rangle\le2
 \qquad\text{for every positive root }\beta,
$$
Theorem \ref{thm:perfect-unique} identifies the two bases.
\end{proof}

Denote the corresponding compact MV cycle by
$Z_{\tau,D}\subseteq\overline{\Gr^{2\varpi_m}}$.

\begin{prop}\label{prop:fusion-cycle}
The special fibre is reduced, and its fundamental cycle is
\begin{equation}\label{eq:fusion-cycle}
 [\mathcal Y_{\lambda,\mu,0}]
 =\sum_{(\tau,D_0)\in\mathcal C(C_\lambda,C_\mu)}[Z_{\tau,D_0}].
\end{equation}
This equality is in the free group of $d$-dimensional cycles on $X$;
the coefficients are the lengths at the generic points.
\end{prop}
\begin{proof}
By Theorem \ref{thm:bd-model}, 
$\mathcal Y_{\lambda,\mu}$ is
exactly the BD fusion of the MV cycles $Y(\lambda)$ and
$Y(\mu)$. 
By \cite[Section 5.2 and Lemma 7.10]{BKK}, 
the multiplicity of $[Z_{\tau,D_0}]$ in $[\mathcal Y_{\lambda,\mu,0}]$ is exactly the multiplicity in \eqref{eq:minortoTL}.

By Proposition \ref{prop:fusion-rational}, the fibre $\mathcal Y_{\lambda,\mu,0}$ is pure of
dimension $d$, these are all its irreducible components.
Finally, the fibre is Cohen--Macaulay by
Proposition \ref{prop:fusion-rational}, hence has no embedded associated primes. Generic reducedness therefore implies
reducedness.
\end{proof}

\begin{theorem}\label{step:(4)}
We have 
$
 \mathcal Y_{\lambda,\mu,0}
 \subseteq\mathcal Y_{\lambda\vee\mu,\lambda\wedge\mu,0}$.
\end{theorem}
\begin{proof}
First, we briefly  explain  that $\mathcal C(C_\lambda,C_\mu)
\subseteq\mathcal C(C_{\lambda\vee\mu},C_{\lambda\wedge\mu})$. 
Write $C_\lambda=\{c_1<\cdots<c_m\}$ and
$C_\mu=\{d_1<\cdots<d_m\}$. Because 
$(\lambda\vee\mu)'_j=\max(\lambda'_j,\mu'_j)$ and
$(\lambda\wedge\mu)'_j=\min(\lambda'_j,\mu'_j)$, we have
$$
C_{\lambda\vee\mu}=\{\min(c_j,d_j):1\le j\le m\},
\qquad
C_{\lambda\wedge\mu}=\{\max(c_j,d_j):1\le j\le m\}.
$$
Then it is readily checked that 
$$C_{\lambda\vee\mu}\cap C_{\lambda\wedge\mu}=C_\lambda\cap C_\mu, \ \ \ \ \ C_{\lambda\vee\mu}\mathbin{\triangle}C_{\lambda\wedge\mu}=C_\lambda\mathbin{\triangle}C_\mu.$$
With the above observations, it is not hard to verify  that if $(\tau,D)\in\mathcal C(C_\lambda,C_\mu)$, then 
$(\tau,D)$   belongs to $\mathcal C(C_{\lambda\vee\mu},C_{\lambda\wedge\mu})$, and so 
$\mathcal C(C_\lambda,C_\mu)
\subseteq\mathcal C(C_{\lambda\vee\mu},C_{\lambda\wedge\mu})$.

Now, Proposition \ref{prop:fusion-cycle} gives the inclusion of the supports.
Both fibres are reduced, so this is a closed-scheme inclusion. 
\end{proof}

By Proposition \ref{step:(1)--(3)} and Theorem \ref{step:(4)},   
Theorem \ref{thm:main} follows from Theorem \ref{thm:criterion}.

\appendix

\section{Rational singularities}\label{sec:AppA}
\def\Hom{\operatorname{Hom}}

In this section, we prove that $\mathcal{Y}_{\lambda,\mu}$ has rational singularities. 
The proof is separated into several  steps. 

\subsection{Background of spherical varieties}Most of the material can be found in Perrin's survey \cite{Perrin}.

Recall that a normal variety with a connected reductive group $G$ action is $G$-\emph{spherical} if a Borel subgroup $B$ has an open orbit. 
Any $G$-spherical variety has finite many $B$-orbits.  

The set of $B$-orbit closures admits a \emph{weak order} generated by $Y\leq P_iY$ for a minimal parabolic $P_i\supset B$. 
A Borel orbit closure
$Y$ is \emph{one-step multiplicity-free} if 
\begin{equation}\label{eq:piPBYtoPY}
\pi: 
 P_i\times_B Y\longrightarrow P_iY,\qquad Y<P_iY,
\end{equation}
has generic degree one. 
A Borel orbit closure
$Y$ is said to be \emph{multiplicity-free} if every $Y'\geq Y$ is one-step multiplicity-free.

\begin{lemma}\label{lem:toconcRS}
Let $X$ be a $G$-spherical variety and $Y$ be a multiplicity-free $B$-orbit closure. 
Assume the following condition for any $Y'\geq Y$ under the weak order with $Y'\neq GY$: 
\begin{itemize}
    \item[(*)] if $Y'$ contains a closed $G$-orbit $Z$, then $Y'$ has a rational singularity at a point of $Z$.
\end{itemize}
Then $Y$ has rational singularities. 
\end{lemma}
\begin{proof}
The remaining argument is a direct  generalization of \cite[Theorem 4.4.7 and Remark 4.4.8]{Perrin}; see also \cite{BrionOrbit}. 
In loc. cit. the condition was that 
$Y'$ does not contain a $G$-orbit. 
From the proof the condition can be relaxed to our statement. Let us sketch the proof for readers' convenience. 

The proof is by downward induction on dimension. 
In characteristic $0$, the $G$-orbit closure $GY$ always has rational singularities. 
Using \eqref{eq:piPBYtoPY}, one can prove $R^{\geq 1}\pi_*\mathcal{O}_{P_i\times_B Y}=0$. 
By induction, it remains to prove $Y$ is normal. 
The locus where $Y$ is not normal is $G$-invariant. Thus its support contains a closed $G$-orbit if it is non-empty. 
But this contradicts our assumption. 
\end{proof}

\begin{remark}
In \cite[Theorem 5]{BrionOrbit}, the theorem is stated without assumption (*). 
The proof actually has a gap; see \cite[Remark 3]{BrionMultfree}. 
\end{remark}

\begin{example}
For $1\leq k,k'\leq n$, 
$\Gr(k,n)\times \Gr(k',n)$ is $GL_{n}$-spherical, and all $B$-orbit closures are multiplicity-free.  
See \cite[Proposition 2.6]{AchingerPerrin} and 
\cite[Corollaries 3.20--3.21]{AchingerPerrin}.
\end{example}

\begin{example}
For $h\geq 1$, 
$GL_{2h}/GL_h\times GL_h$ is $GL_{2h}$-spherical, and all $B$-orbit closures are multiplicity-free. 
See \cite{Wyser16}. Actually, we reduce to above by embedding $GL_{2h}/GL_h\times GL_h$ in $\Gr(h,2h)\times \Gr(h,2h)$. 
\end{example}

We say a variety over $\mathbb{A}^1$ is \emph{horizontal} if its image is dense in $\mathbb{A}^1$. 
We will frequently use the following argument. 

\begin{remark}\label{rmk:frequentarg}
Assume $X$ is a $G$-spherical variety with all $B$-orbit closures multiplicity-free. Recall $\mathbb{U}=\mathbb{A}^1\setminus \{0\}$. 
Assume we have a $G\times \mathbb{G}_m$-variety  $\mathcal{X}$ over $\mathbb{A}^1$ such that 
$$\mathcal{X}_{\mathbb{U}}:=\mathbb{U}\times_{\mathbb{A}^1}\mathcal{X}\cong \mathbb{U}\times X$$
as ($G\times \mathbb{G}_m$)-varieties over $\mathbb{U}$. 
If $\mathcal{X}$ is normal, then $\mathcal{X}$ is ($G\times \mathbb{G}_m$)-spherical since it contains a dense ($B\times \mathbb{G}_m$)-orbit in $\mathcal{X}_{\mathbb{U}}$. 
Moreover, all its horizontal ($B\times\mathbb{G}_m$)-orbit closures are multiplicity-free, as if $Y$ is horizontal, so is every $Y'\geq Y$ under the weak order. 
\end{remark}

\subsection{The family $\mathcal{A}_h$}
Let $h\geq 0$. Define $\mathcal{A}_h$ to be the closure of 
\begin{equation}\label{eq:Ahcirc}
\mathcal{A}_h^\circ
=\{(s,T):s\neq 0,T^2=sT,\operatorname{rank}(T)=\operatorname{rank}(T-sI_{2h})=h\}
\end{equation}
in $\mathbb{A}^1\times \operatorname{End}(\mathbb{C}^{2h})$. 
Note that $\mathcal{A}_h^\circ$ is a single ($GL_{2h}\times \mathbb{G}_m$)-orbit, say, of $\left[\begin{matrix}
I_h & 0 \\ 0 & 0
\end{matrix}\right]$. Thus
\eqref{eq:Ahcirc} is isomorphic to 
$$\mathcal{A}_h^\circ
\cong \mathbb{U}\times (GL_{2h}/GL_h\times GL_h).$$

\begin{prop}\label{prop:descofAh}
We have 
$$\mathcal{A}_h = \{(s,T):T^2=sT,\operatorname{rank}(T)\leq h,\operatorname{rank}(T-sI_{2h})\leq h\}_{\rm red}$$
and the reduced zero fiber is the square-zero matrices
$
\{T\in \operatorname{End}(\mathbb{C}^{2h}):T^2=0\}_{\rm red}$.
\end{prop}
\begin{proof}
When $s\neq 0$, if we are given $T^2=sT$, then 
$$\operatorname{rank}(T)\leq h,
\operatorname{rank}(T-sI_{2h})\leq h\iff 
\operatorname{rank}(T) = \operatorname{rank}(T-sI_{2h})=h,$$
since $\operatorname{im}(T)\oplus \operatorname{im}(T-sI_{2h})=\mathbb{C}^{2h}$. 
When $s=0$, if we are given $T^2=0$, then 
$$\operatorname{rank}(T)\leq h,
\operatorname{rank}(T-sI_{2h})
=\operatorname{rank}(T)\leq h$$
is redundant, since $\operatorname{im}(T)\subseteq \operatorname{ker}(T)$. 
It remains to show any $T$ with $T^2=0$ is contained in $\mathcal{A}_h$. 
By conjugating by an element of $GL_{2h}$, we can assume $T$ is block-diagonal with diagonal blocks 
$
\left[\begin{matrix}
0 & 1 \\ 0 & 0
\end{matrix}\right]$ and 
$\left[\begin{matrix}
0 & 0 \\ 0 & 0
\end{matrix}\right]$. 
Since
$$
\left[\begin{matrix}
s & 1 \\ 0 & 0
\end{matrix}\right]
\stackrel{s\to 0}\longrightarrow
\left[\begin{matrix}
0 & 1 \\ 0 & 0
\end{matrix}\right],\qquad 
\left[\begin{matrix}
s & 0 \\ 0 & 0
\end{matrix}\right] 
\stackrel{s\to 0}\longrightarrow
\left[\begin{matrix}
0 & 0 \\ 0 & 0
\end{matrix}\right],
$$
the proposition follows. 
\end{proof}

\begin{lemma}
The variety $\mathcal{A}_h$ is normal. 
\end{lemma}
\begin{proof}
Firstly, let us write
$\mathcal{A}_h$ as a GIT quotient 
$\mathfrak{A}_h/GL_h$ for 
$$\mathfrak{A}_h=\{(s,u,v)
:vu=sI_h\}
\subset 
\mathbb{A}^1\times \Hom(\mathbb{C}^h,\mathbb{C}^{2h})\times 
\Hom(\mathbb{C}^{2h},\mathbb{C}^{h}).$$
Over $s\neq 0$, the $GL_h$-action is free and we have
$$(\mathbb{U}\times_{\mathbb{A}^1}\mathfrak{A}_h)/GL_h \cong \mathcal{A}_h^\circ.$$
So $\mathbb{U}\times_{\mathbb{A}^1}\mathfrak{A}_h$ has dimension $3h^2+1$ (codimension $h^2$). 
Second, at $s=0$, $(0,u,v)\mapsto (\operatorname{im}(u)\subset \ker(v))$ defines a stratification
$$\{0\}\times_{\mathbb{A}^1}\mathfrak{A}_h
=\bigsqcup_{0\leq a,b\leq h} X_{a,b}$$
where $X_{a,b}$ is an open subset of the total space of 
$\mathcal{S}_1^{\oplus h}\oplus 
\mathcal{Q}_2^{\vee,\oplus h}
$ over partial flag variety $\operatorname{Fl}(h-a,h+b,2h)$ for the first universal subbundle $\mathcal{S}_1$ and the second universal quotient bundle $\mathcal{Q}_2$. 
Thus $\{0\}\times_{\mathbb{A}^1}\mathfrak{A}_h
$ has dimension at most 
$$\begin{aligned}
&\quad  \dim \operatorname{Fl}(h-a,h+b,2h)+\operatorname{rank}(\mathcal{S}_1^{\oplus h}\oplus 
\mathcal{Q}_2^{\vee,\oplus h})\\
& = (h+b)(h-b)+(h-a)(a+b) +h(2h-a-b)\\
& = 3h^2-b^2-a^2-ab\leq 3h^2.    
\end{aligned}
$$
Notice $\mathfrak{A}_h$ is defined by $h^2$ equations, the dimension approximation proves there is no irreducible component that lies in the zero fiber. 
Since $\mathbb{U}\times_{\mathbb{A}^1}\mathfrak{A}_h$ is irreducible, $\mathfrak{A}_h$ is a complete intersection thus Cohen--Macaulay.

Invariant theory tells that $\mathfrak{A}_h/\!/GL_h:=\operatorname{Spec}\mathbb{C}[\mathfrak{A}_h]^{GL_h}$ is generated by $s,uv$. This defines a closed embedding 
$$\mathfrak{A}_h/\!/GL_h\subset \mathbb{A}^1\times \operatorname{End}(\mathbb{C}^{2h}).$$
Since the GIT quotient is isomorphic to $\mathcal{A}_h$ at $s\neq 0$ and the GIT quotient is integral, we can conclude $\mathcal{A}_h=\mathfrak{A}_h/\!/GL_h$.

Using Jacobian criterion, 
$\mathfrak{A}_h$ is smooth over the open locus where $u$ is injective or $v$ is surjective. 
Precisely, the differential 
$$d(vu-sI_{h}) = v\cdot du+dv\cdot u-ds\cdot I$$
is surjective. 
Note that its complement (i.e. $s=0$ and $a,b\geq 1$) has codimension $\geq 4\geq 2$. 
Thus $\mathfrak{A}_h$ is normal using Serre's criterion.
As a GIT quotient, $\mathcal{A}_h$ is also normal. 
\end{proof}

By Remark \ref{rmk:frequentarg}, 
the variety $\mathcal{A}_h$ is ($GL_{2h}\times \mathbb{G}_m$)- spherical.

\begin{theorem}\label{thm:AhhasRS}
Any horizontal $(B_{2h}\times \mathbb{G}_m)$-orbit closure of $\mathcal{A}_h$ has rational singularities. 
\end{theorem}

We will prove it in Section \ref{sec:proofofThRSofA}.

\subsection{Rewrite $\mathcal{X}$}
Recall 
$$\Gr(2k,2n) = \{f\in \Hom(\mathbb{C}^{2n},\mathbb{C}^{2m}): 
f\text{ is surjective}\}/GL_{2m}. $$
We will write $\mathcal{X}$ as a quotient.

\begin{lemma}
We have 
$$
\mathcal{X}=\{(s,T,f)\in \mathcal{A}_m\times \Hom(\mathbb{C}^n,\mathbb{C}^{2m})
: (f,Tf)\text{ is surjective}\}/GL_{2m}. $$
\end{lemma}
\begin{proof}
Let $R_s=\mathbb{C}[t]/(t^2-ts)$.
From the proof of Theorem \ref{thm:bd-model}, we can write
$$\mathcal{X}
=\left\{(s,L): 
\begin{matrix}\text{$L$ is a $\mathbb{C}[t]$-submodule of $R_s^{\oplus n}$ of codimension $2m$}\\
\text{with }\operatorname{rank}(t|_{R_s^{\oplus n}/L})\leq m,\,
\operatorname{rank}((t-s)|_{R_s^{\oplus n}/L})
\leq m
\end{matrix}
\right\}.$$
Equivalently, 
$$\mathcal{X} = \left\{
(s,T,F):
\begin{matrix}
T\in \operatorname{End}(\mathbb{C}^{2m}),
\operatorname{rank}(T)\leq m,
\operatorname{rank}(T-sI_{2m})\leq m\\
\text{$F:R_s^{\oplus n}\to (\mathbb{C}^{2m},T)$
is surjective and $\mathbb{C}[t]$-linear}.
\end{matrix}
\right\}/GL_{2m}. $$
We can naturally identify 
$\mathbb{C}\oplus \mathbb{C}\cong R_s$ by
$(a,b)\mapsto a+bt$. Then any such $F$ can be written as $(f,Tf): \mathbb{C}^{2n}\to \mathbb{C}^{2m}$. 
\end{proof}

Let us identify the zero fiber of $\mathcal{X}$ as 
$$\{V: tV\subseteq V\}\subseteq \Gr(2k,(\mathbb{C}[t]/(t^2))^{\oplus n}).$$
In the following, we will identify a $B^-$-equivariant neighborhood of the torus fixed point 
$$P_{a,b}:=\mathbb{C}^{a}\oplus t\mathbb{C}^b\in \mathcal{X}_0$$
for $a\leq b$ with $a+b=2k$, where we identify $\mathbb{C}^a\subseteq \mathbb{C}^b\subseteq \mathbb{C}^n$. 
Let us consider the  neighborhood $U_{a,b}$ of $P_{a,b}$ to be the locus of $(s,T,f)$ such that  
$$
\{f(e_{i}):a<i\leq n\}\cup 
\{Tf(e_{j}):b<j\leq n\}$$
forms a basis of $\mathbb{C}^{2m}$.
Let us fix such a triple $(s,T,f)$. 
For $a<i\leq b$, we can replace uniquely 
$f(e_i)$ by $f(e'_i)$ for 
$e'_i \in e_i+\operatorname{span}(e_{\imath}:b<\imath\leq n)$
such that $Tf(e_j)$ does not appear in the expansion of $Tf(e'_i)$. 
Since $T^2=sT$, it is easy to see $f(e_i)$ ($b<i\leq n$) does not appear in the expansion of $Tf(e'_i)$ either. 
We partition the standard basis 
$$
\underbrace{e_1,\ldots,e_a}_a,
\underbrace{e_{a+1},\ldots,e_b}_{b-a},
\underbrace{e_{b+1},\ldots,e_n}_{n-b}
\text{ of $\mathbb{C}^n$}
$$
and the new basis 
$$\underbrace{f(e'_{a+1}),\ldots,f(e'_b)}_{b-a},
\underbrace{f(e_{b+1}),\ldots,f(e_n)}_{n-b},
\underbrace{Tf(e_{b+1}),\ldots,Tf(e_n)}_{n-b}
\text{ of $\mathbb{C}^{2m}$}. 
$$
Thus $(s,T,F)$ can be uniquely  represented as 
$$
T = \left[
\begin{matrix}
T' & 0 & 0 \\
0& 0& 0\\
0& I_{n-b}& sI_{n-b}
\end{matrix}\right],
\quad (s,T')\in \mathcal{A}_{(b-a)/2},\qquad 
f = \left[
\begin{matrix}
X & I_{b-a}& 0 \\
Y & Z  & I_{n-b}\\
W & 0& 0
\end{matrix}\right].$$
This proves 
$$U_{a,b}\cong \mathcal{A}_{(b-a)/2}\times (\text{an affine space}). $$

\medbreak
Next, let us identify the $B^-$-action over $U_{a,b}$. 
Let 
$$N=\left[\begin{matrix}
I_a & 0& 0 \\
* & I_{b-a}& 0\\
* & *& I_{n-b}
\end{matrix}\right]\subset B^-_{n}\subset GL_{n}.$$
Note that 
$$
\left[\begin{matrix}
X & I_{b-a}& 0 \\
Y & Z  & I_{n-b}\\
W & 0& 0
\end{matrix}\right]
\left[\begin{matrix}
I_a & 0& 0 \\
a & I_{b-a}& 0\\
b & c& I_{n-b}
\end{matrix}\right]
=
\left[\begin{matrix}
X+a& I_{b-a}& 0 \\
Y+Za+b& Z+c  & I_{n-b}\\
W & 0& 0
\end{matrix}\right].
$$
This proves the neighborhood $U_{a,b}$ is $N$-equivariant. 
Moreover, we can cancel $X,Y,Z$ using a unique choice of $a,b,c$. This proves 
$$U_{a,b}\cong \mathcal{A}_{(b-a)/2}\times N\times \mathbb{M}_{(n-b)\times a}. $$
Next, note that 
$$
\left[\begin{matrix}
0 & I_{b-a}& 0 \\
0 & 0  & I_{n-b}\\
W & 0& 0
\end{matrix}\right]
\left[\begin{matrix}
x & 0& 0 \\
0 & y& 0\\
0 & 0& z
\end{matrix}\right]
=
\left[\begin{matrix}
0& y& 0 \\
0& 0  & z\\
Wx & 0& 0
\end{matrix}\right]
=
\left[\begin{matrix}
y & 0& 0 \\
0 & z &0 \\
0 & 0& z
\end{matrix}\right]
\left[\begin{matrix}
0 & I_{b-a}& 0 \\
0 & 0  & I_{n-b}\\
z^{-1}Wx & 0& 0
\end{matrix}\right]
$$
and 
$$
\left[
\begin{matrix}
y & 0 & 0 \\
0& z& 0\\
0& 0 & z
\end{matrix}\right]
\left[
\begin{matrix}
T' & 0 & 0 \\
0& 0& 0\\
0& I_{n-b}& sI_{n-b}
\end{matrix}\right]
\left[
\begin{matrix}
y & 0 & 0 \\
0& z& 0\\
0& 0 & z
\end{matrix}\right]^{-1}
=
\left[
\begin{matrix}
yT'y^{-1} & 0 & 0 \\
0& 0& 0\\
0& I_{n-b}& sI_{n-b}
\end{matrix}\right].
$$
Thus $U_{a,b}$ is not only invariant under $N$, but also invariant under $B^-\subset GL_{n}$. 
Moreover, any ($B^-\times \mathbb{G}_m$)-orbit (closure) of $U_{a,b}$ is of the form
$$N\times A\times M$$
for a ($B^-_{b-a}\times\mathbb{G}_m$)-orbit (closure) of $A\subset \mathcal{A}_{(b-a)/2}$ and a ($B^-_{n-b}\times B^-_a$)-orbit (closure) of  $M\subseteq \mathbb{M}_{(n-b)\times a}$.

\begin{theorem}\label{prop:RSofYlambdaYmu}
The variety $\mathcal{Y}_{\lambda,\mu}$ has rational singularities. 
\end{theorem}
\begin{proof}
The fiber of $\mathcal{X}=\overline{\aGr^{\varpi_m,\varpi_m}_{\mathbb{A}}}$ is either
$\overline{\aGr^{\varpi_m}}\times \overline{\aGr^{\varpi_m}}$ or
$\overline{\aGr^{2\varpi_m}}$, 
which has rational singularities since they are Schubert varieties \cite{Kum}. 
The morphism $\mathcal{X}\to\mathbb A^1$ is flat.
Its base is smooth, and all its fibres have rational
singularities. Hence $\mathcal X$ has rational singularities
by \cite[Th\'eor\`eme V]{Elkik}.
In particular, $\mathcal{X}$ is normal. 
By Remark \ref{rmk:frequentarg}, 
$\mathcal{X}$ is spherical. 
We prove a stronger result that any horizontal ($B^-\times \mathbb{G}_m$)-orbit closure has rational singularities. 
We need to check (*) in Lemma \ref{lem:toconcRS}

Over $\mathcal{X}$, the closed ($GL_n\times \mathbb{G}_m$)-orbits are all in the zero fiber, and thus projective. 
In particular, it contains a torus fixed point. Let us pick the point to be the unique $B\subset GL_n$-invariant fixed point, i.e. $P_{a,b}$. 
Let $Z$ be a horizontal ($B^-\times \mathbb{G}_m$)-orbit closure of $\mathcal{X}$. 
By above discussion 
$$Z\cap U_{a,b}\cong N\times A\times M$$
where $A$ is a horizontal ($B_{b-a}\times \mathbb{G}_m$)-orbit closure and $M$ is a ($B^-_{n-b}\times B^-_a$)-orbit closure, i.e. a matrix Schubert variety. 
By Theorem \ref{thm:AhhasRS} and \cite[Theorem 2.4.3 and paragraph after it]{KM}, each factor has rational singularities.
\end{proof}

\subsection{Proof of Theorem \ref{thm:AhhasRS}}
\label{sec:proofofThRSofA}

Let us consider the open loci 
$$\mathcal{A}_h^{\prime}
=\{(s,T)\in \mathcal{A}_h: \operatorname{rank}(T)=h\}. $$
Note that its zero fiber contains only one closed $G$-orbit, i.e. the nilpotent orbit of Jordan type $(2^h)$. 
It is easy to compute that the dimension of the nilpotent orbit is $2h^2$. 

\begin{lemma}
Any horizontal ($B_{2h}^-\times \mathbb{G}_m$)-orbit closure $Y$ of $\mathcal{A}_h^{\prime}$ has rational singularities. 
\end{lemma}
\begin{proof}
Clearly $\mathcal{A}'_h$ is a spherical variety. 
This follows immediately from Lemma \ref{lem:toconcRS}. 
Actually, any horizontal $Y'\neq GY$ cannot contain the zero fiber by dimension reason. 
\end{proof}

Let $e\in \operatorname{End}(\mathbb{C}^{2h})$ with 
$e^2=e$, and $\operatorname{rank}(e)=h$. 
We define 
$$O_e = \{(s,sbeb^{-1}):s\in \mathbb{C}^\times,b\in B_{2h}\}
\subset \mathcal{A}_h.$$
Any horizontal ($B_{2h}\times\mathbb{G}_m$)-orbit closure of $\mathcal{A}_h$ is of the form $\overline{O_e}$. 
Assume $\overline{O_f}=\mathcal{A}_h$ for some $f$.  

Let us construct 
$$\tilde{e} = \left[\begin{matrix}
e & 0 \\
0 & f
\end{matrix}\right]\in \operatorname{End}(\mathbb{C}^{4h}).$$

\begin{lemma}
There is an open embedding 
$$\overline{O_e}\times GL_{2h}\longrightarrow \overline{O_{\tilde{e}}}\cap \mathcal{A}_{2h}^{\prime},\qquad 
(s,T)\times g\longmapsto 
\left(s,
\left[\begin{matrix}
    T & g\\ 0 & g^{-1}(sI_{2h}-T)g
\end{matrix}\right]
\right).$$
\end{lemma}
\begin{proof}
Let us first check the image lies in $\overline{O_{\tilde{e}}}$. 
Let us consider the open locus 
$$\Omega = \{(s,T)\times g:
(s,T)\in O_e, 
(s,g^{-1}(sI_{2h}-T)g)\in O_{f}
\}\subset \overline{O_e}\times GL_{2h}.$$
This is nonempty since $G\cdot O_e=G\cdot O_f=\mathbb{U}\times_{\mathbb{A}^1}\mathcal{A}_h$. 
Let $T'=g^{-1}(sI_{2h}-T)g$. 
We have 
$$
\left[\begin{matrix}
T & g \\
0 & T' 
\end{matrix}\right]
=\left[\begin{matrix}
I_{2h} & W \\
0 & I_{2h} 
\end{matrix}\right]
\left[\begin{matrix}
T & 0 \\
0 & T' 
\end{matrix}\right]
\left[\begin{matrix}
I_{2h} & -W \\
0 & I_{2h} 
\end{matrix}\right],\qquad 
W = \frac{1}{s^2}(gT'-Tg).$$
In the above identity, one might use 
$$T^2=sT,\qquad (T')^{2} = sT',\qquad Tg+gT'=sg.$$
This proves the image of $\Omega$, thus the image of $\overline{O_e}\times GL_{2h}$ is contained in $\overline{O_{\tilde{e}}}$. 

Next, note that any 
$(s,\tilde{T})\in \overline{O_{\tilde{e}}}$, we can write 
$$\tilde{T}=\left[\begin{matrix}
T & g \\ 0 & T'
\end{matrix}\right].$$
Actually, this is true for $(s,\tilde{T})\in O_{\tilde{e}}$ and 
the vanishing of lower-left block is a closed condition. 
Let us restrict to the open loci $\Sigma$ of $\overline{O_{\tilde{e}}}$ when $g$ is invertible. 
Note that $\operatorname{rank}(\tilde{T})\geq 2h$, but $\operatorname{rank}(\tilde{T})\leq 2h$ by Proposition \ref{prop:descofAh}. 
Thus $\Sigma\subset\mathcal{A}_{2h}'$. 
The equation $\tilde{T}^2= s\tilde{T}$ gives 
$$T'=g^{-1}(sI_{2h}-T)g. $$
This proves $\overline{O_e}\times GL_{2h}\cong \Sigma\subseteq \overline{O_{\tilde{e}}}\cap \mathcal{A}_{2h}^{\prime}$. 
\end{proof}

By the two lemmas above, $\overline{O_e}\times GL_{2h}$ has rational singularities. 
Since $GL_{2h}$ is smooth, $\overline{O_e}$ has rational singularities. 
Now Theorem \ref{thm:AhhasRS} follows.

\section{Uniqueness of perfect bases}\label{sec:perfect}
We prove the uniqueness statement used in Lemma \ref{lem:tlmv}.

Let $L(\lambda)$ be an irreducible finite-dimensional representation
of a connected complex reductive group. Fix a highest vector
$v_\lambda$ and Chevalley operators $e_i,f_i$, with
$h_i=[e_i,f_i]$; $i$ ranges over the simple roots $\alpha_i$.
For $0\ne v\in L(\lambda)$, put
$$
 \varepsilon_i(v)=\max\{a\ge0:e_i^av\ne0\}.
$$

\begin{definition}\label{def:perfect}
A weight basis $\mathcal B$ containing $v_\lambda$ of $L(\lambda)$ is \emph{perfect}
if, for each $i$, there is a map
$\widetilde e_i:\mathcal B\to\mathcal B\sqcup\{0\}$ such that
\begin{enumerate}
\item $\widetilde e_i b=0$ if and only if $e_i b=0$;
\item $\widetilde e_i$ is injective on its nonzero domain;
\item if $a=\varepsilon_i(b)>0$, then
\begin{equation}\label{eq:perfect-definition}
 e_i b = \varepsilon_i(b)\widetilde e_i b+(\cdots),\qquad (\cdots)
\in\ker e_i^{a-1}.
\end{equation}
\end{enumerate}
\end{definition}

\begin{theorem}\label{thm:perfect-unique}
If the dominant weight $\lambda$ satisfies
$$
 \langle\lambda,\beta^\vee\rangle\le2
 \qquad\text{for every positive root }\beta,
$$
then $L(\lambda)$ has at most one 
perfect basis.
Here $\beta^\vee$ denotes the coroot.
\end{theorem}
\begin{proof}
For two perfect bases $\mathcal B,\mathcal B'$, the comparison
theorem \cite[Corollary 5.39]{BK} gives a weight-preserving
bijection
$$
 \psi:\mathcal B\longrightarrow\mathcal B',
 \qquad
 \psi(\widetilde e_i b)=\widetilde e_i\psi(b),
 \qquad \psi(0)=0.
$$
Note that the condition is equivalent to saying 
$$
 -2\le\langle\mu,h_i\rangle\le2
 \qquad\text{for every weight }\mu\text{ and every }i.
$$
Now let us induct downwards from the highest weight. Let $b\in\mathcal B$
have weight $\mu\ne\lambda$, and suppose the assertion holds
at all higher weights.

If $\langle\mu,h_i\rangle=-2$ for some $i$, then
$\mathfrak{sl}_2$-theory gives an injection
$$
 e_i^{(2)}:=\frac{e_i^2}{2}:L(\lambda)_\mu\longrightarrow
                 L(\lambda)_{\mu+2\alpha_i}.
$$
Both $b$ and $\psi(b)$ have $\varepsilon_i=2$. 
We have 
$$
 e_i^{(2)}b=\widetilde e_i^2b
 =\psi(\widetilde e_i^2b)=e_i^{(2)}\psi(b).
$$
Injectivity gives $b=\psi(b)$.

Otherwise $\langle\mu,h_i\rangle\ge-1$ for every $i$.
The weight bound gives
$$
 e_i^2L(\lambda)_\mu=0,
 \qquad \varepsilon_i(b)\le1
 \qquad\text{for every }i.
$$
Since $\varepsilon_i(b)\leq 1$, we have 
$e_ib = \varepsilon_i(b){\widetilde e_i}b$ 
by \eqref{eq:perfect-definition}. 
Thus
$$
 e_i b=\widetilde e_i b
 =\psi(\widetilde e_i b)
 =\widetilde e_i\psi(b)=e_i\psi(b)
 \qquad\text{for every }i.
$$
It follows that
$$
 b-\psi(b)\in L(\lambda)_\mu\cap\bigcap_i\ker e_i
 =L(\lambda)_\mu\cap \mathbb{C}v_\lambda=0.
\qedhere$$
\end{proof}

\end{document}